\documentclass[smallextended]{svjour3}
\usepackage{graphicx}
\usepackage{fix-cm}
\usepackage{amssymb}
\usepackage{amsmath}
\usepackage{comment}
\usepackage{bm}
\usepackage{multirow}
\usepackage[numbers]{natbib}
\usepackage{microtype}

\numberwithin{equation}{section}
\allowdisplaybreaks

\newtheorem{thm}{Theorem}

\begin{document}

\title{On a fractional binomial process 
}


\author{Dexter O. Cahoy         \and
        Federico Polito 
}


\institute{D.O. Cahoy \at
Department of Mathematics and Statistics, College of Engineering and Science \\
	Louisiana Tech University, Ruston, LA, USA \\
  \email{dcahoy@latech.edu}
      \and
           F. Polito \at
      Dipartimento di Matematica \\
	Universit\`a degli studi di Roma ``Tor Vergata'', Rome, Italy \\
\email{polito@nestor.uniroma2.it}
}

\date{Received: date / Accepted: date}

\maketitle

\begin{abstract}
The classical binomial process has been studied by \citet{jakeman}
		(and the references therein) and has been used to characterize a
		series of radiation states in quantum optics. In particular, he studied a classical birth-death process
		where the chance of birth is proportional to the difference between a larger fixed number  and the number of
		individuals present. It is shown that at large times, an equilibrium is reached which follows a binomial process.
		In this paper, the classical binomial process is generalized using the techniques of fractional calculus and is called
		the fractional binomial process. The fractional binomial process is shown to preserve the binomial limit at large times
		while expanding the class of models that include non-binomial fluctuations (non-Markovian) at regular and small  times.  As a
		direct consequence, the generality of the fractional binomial model makes the proposed model more  desirable than its classical counterpart in describing real physical
		processes. More statistical properties are also derived.
		
\end{abstract}

\section{Introduction}
	The classical binomial process  has been studied by  \citet{jakeman} and has been used to model fluctuations in a train of events in quantum optics. 
	Recall that the classical binomial process $ \mathcal{N}(t)$, $t\ge 0 $, with birth rate $\lambda>0$ and
		death rate $\mu>0$, has  state probabilities
		$p_n(t) = \Pr \{ \mathcal{N}(t) = n | \mathcal{N}(0) = M \}$ which solve the following Cauchy problem: 
		\begin{align}
			\label{state}
			\begin{cases}
				\frac{\mathrm d}{\mathrm dt} p_n(t) = \mu (n+1) p_{n+1}(t) - \mu n p_n(t) - \lambda
				(N-n) p_n(t) \\
				\qquad \qquad + \lambda (N-n+1) p_{n-1}(t), \qquad \qquad \qquad 0 \leq n \leq N, \\
				p_n(0) =
				\begin{cases}
					1, & n=M, \\
					0, & n \neq M.
				\end{cases}
			\end{cases}
		\end{align}
		The initial number of individuals is $M \geq 1$, and $N \geq M$.

		Notice that the binomial process has a completely different behaviour compared to the
		classical linear birth-death process. Here the birth rate is proportional to the
		difference between a larger fixed number  and the number of	individuals present while the
		death rate remains linear. The whole evolution of the binomial process develops
		in the region $[0, N]$.
		Furthermore it is shown that at large times, an equilibrium is reached and displays a binomial distribution.
	
		From \eqref{state}, it is straightforward to realise that the generating function
		\begin{align}
			Q(u,t) = \sum_{n=0}^N (1-u)^n p_n(t), \qquad |1-u| \leq 1,
		\end{align}
		is the solution to
		\begin{align}
			\label{fgp}
			\begin{cases}
				\frac{\partial}{\partial t} Q(u,t) = -\mu u \frac{\partial}{\partial u} Q(u,t) - \lambda u
				(1-u) \frac{\partial}{\partial u} Q(u,t) -\lambda N u Q(u,t), \\
				Q(u,0) = (1-u)^M.
			\end{cases}
		\end{align}
Moreover, \citet{jakeman}  showed that at large times, the evolving population follows a  binomial distribution with parameter $\lambda / (\lambda  + \mu)$.  

In this paper, we  propose a fractional generalisation of the classical binomial process. The fractional generalization includes non-markovian and rapidly dissipating or bursting birth-death processes at small and regular times. We also derive more statistical and related properties of the newly developed fractional stochastic process, which are deemed useful in real applications.  Note that the theory and results presented here may have applications beyond quantum optics and may be of interest in other disciplines. 	As in the preceding works on fractional Poisson process (e.g.\ \citet{laskin}) and other
		fractional point processes (see e.g.\ \citet{cah,pol}), fractionality is obtained by replacing the
		integer-order derivative in the governing differential equations with a fractional-order derivative. 	In particular, we use the Caputo fractional derivative of a well-behaved function $f(t)$ and is defined as
		\begin{align}
			\label{caputo}
			\frac{\mathrm d^\nu}{\mathrm dt^\nu} f(t) = \frac{1}{\Gamma(m-\nu)} \int_0^t \frac{\frac{\mathrm d^m}{
			\mathrm d \tau^m}f(\tau)}{(t-\tau)^{\nu-m+1}}
			\mathrm d\tau, \qquad m=\lceil \nu \rceil,
		\end{align}
	where ``$\lceil y \rceil$" is the smallest integer that is not less than  $y$. Note that the Caputo fractional derivative operator is in practice a convolution of the standard derivative
with a power law kernel  which  adds more memory in the process. This characteristic  is  certainly an improvement  from a physical viewpoint. By simple substitution, we obtain the following initial value problems for the probability generating function and the state
		probabilities:
		\begin{align}
			\label{fgpfrac}
			\begin{cases}
				\frac{\partial^\nu}{\partial t^\nu} Q^\nu(u,t) = -\mu u \frac{\partial}{\partial u} Q^\nu(u,t) - \lambda u
				(1-u) \frac{\partial}{\partial u} Q^\nu(u,t) -\lambda N u Q^\nu(u,t), \\
				Q^\nu(u,0) = (1-u)^M, \qquad \qquad \qquad |1-u|\leq 1,
			\end{cases}
		\end{align}		
		\begin{align}
			\label{statefrac}
			\begin{cases}
				\frac{\mathrm d^\nu}{\mathrm dt^\nu} p_n^\nu(t) = \mu (n+1) p_{n+1}^\nu(t) - \mu n p_n^\nu(t) - \lambda
				(N-n) p_n^\nu(t) \\
				\qquad \qquad \qquad + \lambda (N-n+1) p_{n-1}^\nu(t), & 0 \leq n \leq N, \\
				p_n^\nu(0) =
				\begin{cases}
					1, & n=M, \\
					0, & n \neq M,
				\end{cases}
			\end{cases}
		\end{align}
		where $\nu \in (0,1]$.

We organized the rest of the paper as follows. In Section 2, the statistical properties of the fractional binomial process are derived by solving the preceding initial-value problem.  Section 3 explored the sub-models that are directly extractable from the fractional binomial process. We then conclude the paper  by providing more discussions and future extensions of the study in Section 4.

	\section{Main properties of the fractional binomial process}
		\label{se}
		
		Firstly, we prove a subordination relation which is of fundamental importance to deriving many of our results.

		\begin{thm}
			\label{sub}
			The fractional binomial process $\mathcal{N}^\nu(t)$ has the following one-dimensional representation:
			\begin{align}
				\mathcal{N}^\nu(t) \overset{\text{d}}{=} \mathcal{N}(V_t^\nu), 
			\end{align}
			where $\mathcal{N}(t)$ is a classical binomial process, $V_t^\nu$, $t\ge 0$, is the inverse process
			of the $\nu$-stable subordinator (see e.g.\ \citet{meer}),  $t \ge 0$, and $\nu \in (0,1]$.
			
			\begin{proof}
				Let $\text{Pr} \{ V_t^\nu \in \mathrm ds \} = h(s,t) \, \mathrm ds$
				be the law of the inverse $\nu$-stable subordinator. We now show that
				\begin{align}
					Q^\nu(u,t) = \sum_{n=0}^N (1-u)^n p_n^\nu(t) = \int_0^\infty Q(u,s) \, h(s,t) \, \mathrm ds
				\end{align}
				satisfy the fractional differential equation \eqref{fgpfrac}. We can then write
				\begin{align}
					\frac{\partial^\nu}{\partial t^\nu} \int_0^\infty Q(u,s) h(s,t) \mathrm ds
					= \int_0^\infty Q(u,s) \frac{\partial^\nu}{\partial t^\nu} h(s,t) \mathrm ds.
				\end{align}
				Since it can be easily verified that $h(s,t)$ is a solution to the fractional equation
				\begin{align}
					\frac{\partial^\nu}{\partial t^\nu} h(s,t) = - \frac{\partial}{\partial s} h(s,t),
				\end{align}
				we readily obtain
				\begin{align}
					& \frac{\partial^\nu}{\partial t^\nu} Q^\nu(u,t) \\
					& = - \int_0^\infty Q(u,s)
					\frac{\partial}{\partial s} h(s,t) \mathrm ds \notag \\
					& = \left. - h(s,t) Q(u,s) \right|_{s=0}^\infty + \int_0^\infty h(s,t) \frac{\partial}{\partial s}
					Q(u,s) \mathrm ds \notag \\
					& = \int_0^\infty \left[ -\mu u \frac{\partial}{\partial u} Q(u,s)
					-\lambda u(1-u) \frac{\partial}{\partial u} Q(u,s)
					-\lambda NuQ(u,s) \right] h(s,t) \mathrm ds \notag \\
					& = -\mu u \frac{\partial}{\partial u} Q^\nu(u,t)
					-\lambda u(1-u) \frac{\partial}{\partial u} Q^\nu(u,t)
					-\lambda NuQ^\nu(u,t). \notag
				\end{align}
			\end{proof}
			\begin{flushright} \qed \end{flushright}
		\end{thm}
		
\begin{figure}[h!t!tb!p!]
			\centering
			\includegraphics[height=2.5in, width=2.2in]{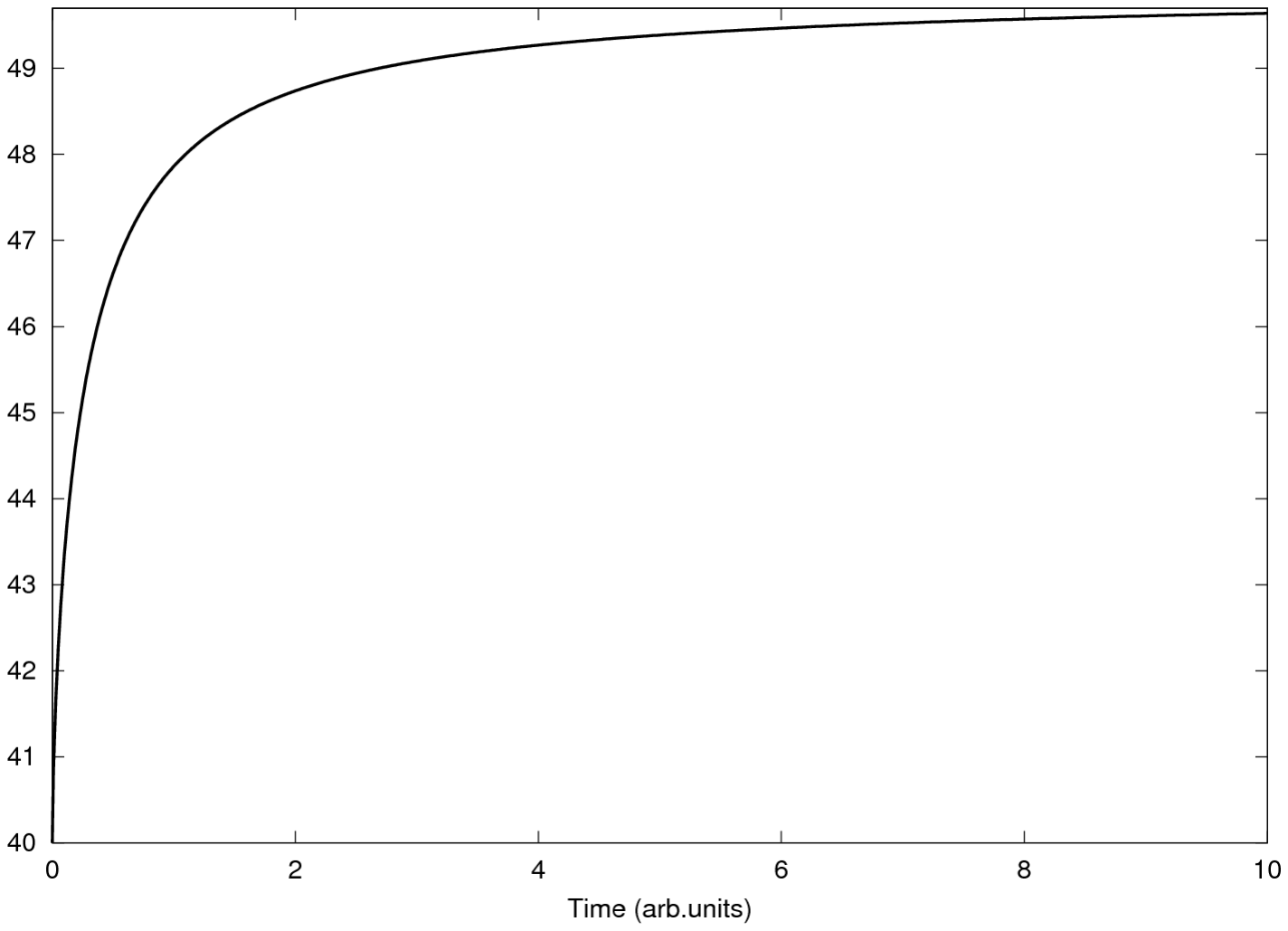}
			\includegraphics[height=2.5in, width=2.2in]{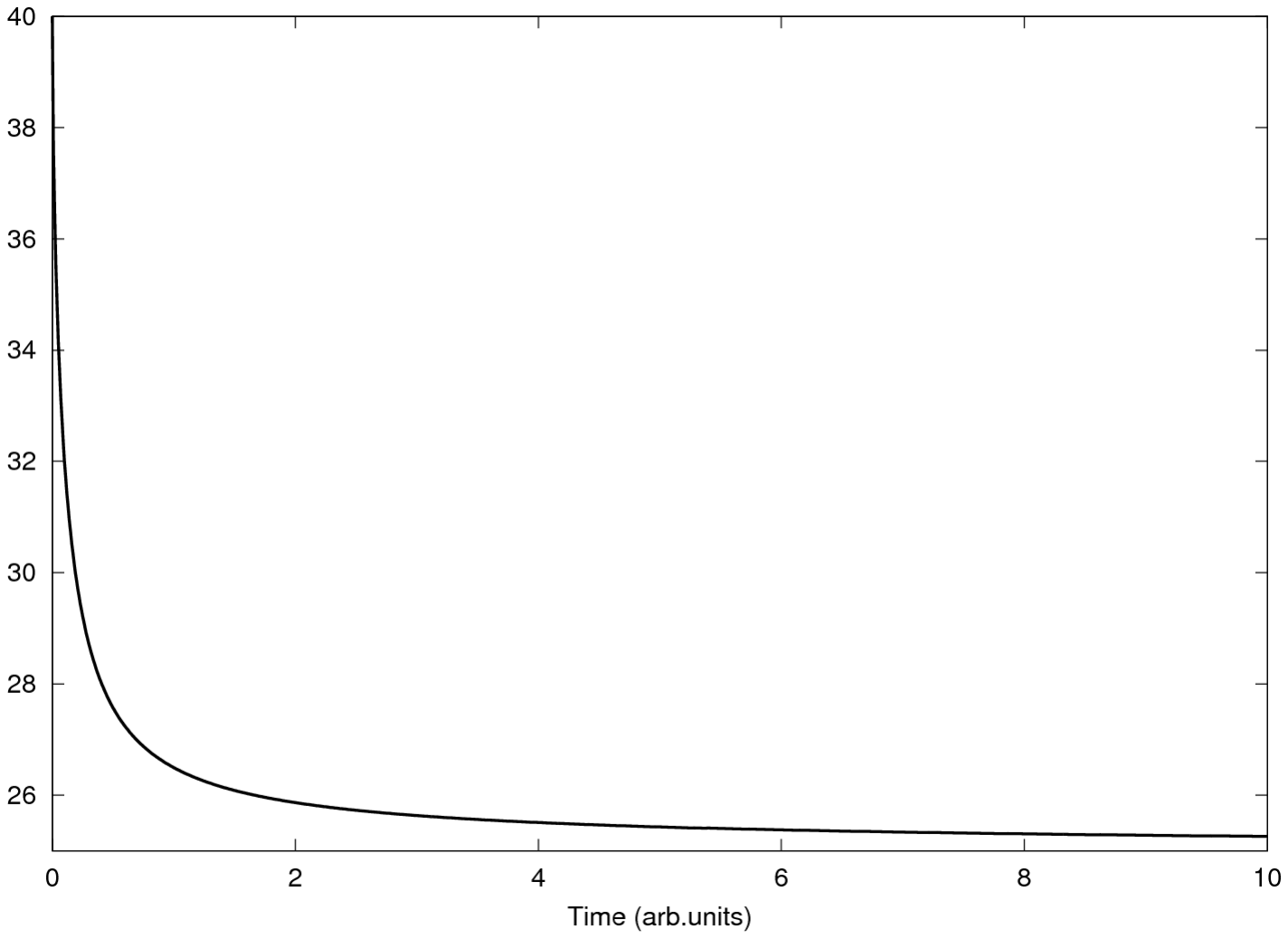}
			\caption{\label{afig}The mean value of the fractional binomial process $\mathbb{E}\,
				\mathcal{N}^\nu(t)$. For both graphs we have $N=100$, $M=40$, $\nu=0.7$. The rates are respectively
				$(\lambda,\mu)=(1,1)$ (left) and $(\lambda,\mu)=(1,3)$ (right).}
		\end{figure}

		In the following theorem, we derive the expected number of individuals $\mathbb{E}\,\mathcal{N}^\nu(t)$ or the expected population size of the fractional binomial process  at any time $t \ge 0$.
		\begin{thm}
			For the fractional binomial process $\mathcal{N}^\nu(t)$, $t \ge 0$, $\nu \in (0,1]$, we have
			\begin{align}
				\label{mea}
				\mathbb{E}\,\mathcal{N}^\nu(t) = \left( M-N\frac{\lambda}{\lambda+\mu} \right)
				E_{\nu,1} \left( -(\lambda+\mu) t^\nu \right) + N\frac{\lambda}{\lambda+\mu},
			\end{align}
	where
	\[
	E_{\alpha, \beta} \left( \xi \right) = \sum\limits_{r=0}^\infty \frac{\xi^r}{\Gamma (\alpha r + \beta) }	
	\]
is the Mittag-Leffler function.			
			\begin{proof}
				By considering that
				\begin{align}
					\label{sharp}
					\left. -\frac{\partial}{\partial u} Q^\nu(u,t) \right|_{u=0} = \mathbb{E}\,\mathcal{N}^\nu(t)
				\end{align}
				and on the base of \eqref{fgpfrac},
				we can write
				\begin{align}
					-\frac{\partial^\nu}{\partial t^\nu} \frac{\partial}{\partial u} Q^\nu(u,t) = {} &
					\mu \left( \frac{\partial}{\partial u} Q^\nu(u,t) + u \frac{\partial^2}{\partial u^2}
					Q^\nu(s,t) \right) \\
					& + \lambda \left( \frac{\partial}{\partial u} Q^\nu(u,t)
					+ u \frac{\partial^2}{\partial u^2} Q^\nu(u,t) \right) \notag \\
					& - \lambda \left( 2 u \frac{\partial}{\partial u} Q^\nu(u,t) + u^2
					\frac{\partial^2}{\partial u^2}Q^\nu(u,t) \right) \notag \\
					& + \lambda N \left( Q^\nu(u,t)
					+ u \frac{\partial}{\partial u} Q^\nu(u,t) \right), \notag
				\end{align}
				thus leading to the Cauchy problem
				\begin{align}
					\label{bel}
					\begin{cases}
						\frac{\mathrm d^\nu}{\mathrm dt^\nu} \mathbb{E} \mathcal{N}^\nu(t) =
						-(\mu +\lambda) \mathbb{E} \mathcal{N}^\nu(t) + \lambda N, \\
						\mathbb{E} \mathcal{N}^\nu(0) = M.
					\end{cases}
				\end{align}
				The solution to \eqref{bel} can be written as (using formula (4.1.65) of \citet{kilbas})
				\begin{align}
					\mathbb{E}\,\mathcal{N}^\nu(t)
					= {} & M E_{\nu,1} \left(-(\lambda+\mu)t^\nu\right) \\
					& + \int_0^t (t-s)^{\nu-1}
					E_{\nu,\nu} \left(-(\lambda+\mu)(t-s)^\nu\right) \lambda N \mathrm ds \notag \\
					= {} & M E_{\nu,1} \left(-(\lambda+\mu)t^\nu\right) + \lambda N \int_0^t y^{\nu-1} E_{\nu,\nu}
					\left(-(\lambda+\mu) y^\nu\right) \mathrm dy \notag \\
					= {} & M E_{\nu,1} \left(-(\lambda+\mu)t^\nu\right) + \lambda N \biggl| -\frac{
					E_{\nu,1}\left( -(\lambda+\mu)y^\nu \right)}{(\lambda+\mu)} \biggr|_0^t \notag \\
					= {} & M E_{\nu,1} \left(-(\lambda+\mu)t^\nu\right) - \frac{\lambda}{\lambda+\mu}
					N \left( E_{\nu,1} \left( -(\lambda+\mu) t^\nu \right) -1 \right) \notag \\
					= {} & \left( M-N\frac{\lambda}{\lambda+\mu} \right)E_{\nu,1}
					\left( -(\lambda+\mu) t^\nu \right) + N\frac{\lambda}{\lambda+\mu}. \notag
				\end{align}
			\end{proof}
			\begin{flushright} \qed \end{flushright}
		\end{thm}
		Figure \ref{afig} shows the mean value \eqref{mea} in both cases $\left[ M-N\lambda/(\lambda+\mu)
		\right] < 0$ and $\left[ M-N\lambda/(\lambda+\mu)
		\right] > 0$ for specific values of the remaining parameters.
		Note also that when $M=N\lambda/(\lambda+\mu)$ the mean value $\mathbb{E}\,
		\mathcal{N}^\nu(t) = N\lambda/(\lambda+\mu)$ is constant.
	
		We now proceed to deriving the variance $\mathbb{V}\text{ar} \, \mathcal{N}^\nu(t)$ of the fractional
		binomial process, starting from the second factorial moment.
		
		\begin{thm}		
					For the fractional binomial process $\mathcal{N}^\nu(t)$, $t \ge 0$, $\nu \in (0,1]$, we have
			\begin{align}
				\label{va}
				\mathbb{V}\text{ar} & \, \mathcal{N}^\nu(t) \\			
				= {} & \left( \frac{\lambda^2 N(N-1)}{(\lambda +\mu)^2} -\frac{2\lambda M(N-1)}{\lambda+\mu}
				+ M(M-1) \right) E_{\nu,1} (-2(\lambda+\mu)t^\nu) \notag \\
				& + \left( \frac{2\lambda^2 N}{(\lambda+\mu)^2} - \frac{\lambda}{\lambda+\mu}
				(N+2M) +M \right)
				E_{\nu,1} (-(\lambda+\mu)t^\nu) \notag \\
				& - \left( M-N\frac{\lambda}{\lambda+\mu} \right)^2 \left( E_{\nu,1}(-(\lambda+\mu)t^\nu) \right)^2
				+ \frac{N \lambda\mu}{(\lambda+\mu)^2}. \notag
			\end{align}
			
			\begin{proof}
				From \eqref{fgpfrac}, we have
				\begin{align}
					\frac{\partial^\nu}{\partial t^\nu} \frac{\partial^2}{\partial u^2} Q^\nu(u,t) = {} &
					-\mu \frac{\partial^2}{\partial u^2} Q^\nu(u,t) -\mu \left( \frac{\partial^2}{\partial u^2}
					Q^\nu(u,t) + u \frac{\partial^3}{\partial u^3} Q^\nu(u,t) \right) \\
					& -\lambda \left( (1-2u) \frac{\partial^2}{\partial u^2} Q^\nu(u,t) -2 \frac{\partial}{\partial u}
					Q^\nu(u,t) \right) \notag \\
					& - \lambda \left( (1-2u) \frac{\partial^2}{\partial u^2}Q^\nu(u,t)
					+(u-u^2)\frac{\partial^3}{\partial u^3} Q^\nu(u,t) \right) \notag \\
					& -\lambda N \frac{\partial}{\partial u} Q^\nu(u,t) - \lambda N \left( \frac{\partial}{\partial u}
					Q^\nu(u,t) + u \frac{\partial^2}{\partial u^2} Q^\nu(u,t) \right). \notag
				\end{align}
				Recalling \eqref{sharp} and the equality
				\begin{align}
					\left. \frac{\partial^2}{\partial u^2} Q^\nu(u,t) \right|_{u=0} = \mathbb{E}
					\left( \mathcal{N}^\nu(t)(\mathcal{N}^\nu(t) -1) \right) = H^\nu(t),
				\end{align}
				we obtain
					\begin{align}
					\label{sarp}
					\frac{\mathrm d^\nu}{\mathrm dt^\nu} H^\nu(t) & =
					-2 \mu H^\nu(t) -2\lambda H^\nu(t) -2\lambda \mathbb{E}\,\mathcal{N}^\nu(t) + 2\lambda N
					\mathbb{E}\,\mathcal{N}^\nu(t) \\
					& = -2 (\lambda+\mu) H^\nu(t) + 2 \lambda (N-1) \mathbb{E}\,\mathcal{N}^\nu(t). \notag
				\end{align}
				By substituting \eqref{mea} into \eqref{sarp}, we arrive at the Cauchy problem
				\begin{align}
					\label{cacio}
					\begin{cases}
						\frac{\mathrm d^\nu}{\mathrm dt^\nu} H^\nu(t) = -2 (\lambda+\mu) H^\nu(t)
						+ 2 \lambda (N-1) \left( M-N\frac{\lambda}{\lambda+\mu} \right)
						E_{\nu,1} \left( -(\lambda+\mu)t^\nu \right) \\
						\qquad \qquad \qquad + 2 \lambda^2 N (N-1) \frac{1}{\lambda+\mu} \\
						H^\nu(0) = M(M-1),
					\end{cases}
				\end{align}
				that can be solved using the Laplace transform $\widetilde{H}^\nu(z) = \int_0^\infty
				e^{-zt} H^\nu(t)\, \mathrm dt$ as follows:
				\begin{align}
					z^\nu \widetilde{H}^\nu(z) &- z^{\nu-1} M(M-1) \\
					= {} & -2 (\lambda+\mu) \widetilde{H}^\nu(z) + 2 \lambda (N-1)
					\left( M-N\frac{\lambda}{\lambda+\mu} \right) \frac{z^{\nu-1}}{z^\nu + (\lambda+\mu)} \notag \\
					& + \frac{1}{z} 2 \lambda^2 N(N-1)\frac{1}{\lambda+\mu}. \notag
				\end{align}
				The Laplace transform then reads
				\begin{align}
					\label{elastic}
					\widetilde{H}^\nu(z) = {} & M(M-1) \frac{z^{\nu-1}}{z^\nu+2(\lambda+\mu)} \\
					& +2\lambda (N-1)
					\left( M-N\frac{\lambda}{\lambda+\mu} \right) \frac{z^{\nu-1}}{(z^\nu+(\lambda+\mu))
					(z^\nu+2(\lambda+\mu))} \notag \\
					& + 2 \lambda^2 N(N-1) \frac{1}{\lambda+\mu} \cdot \frac{z^{-1}}{z^\nu
					+2(\lambda+\mu)} \notag \\
					= {} & M(M-1) \frac{z^{\nu-1}}{z^\nu+2(\lambda+\mu)} \notag \\
					& + \frac{2\lambda(N-1)}{\lambda+\mu}
					\left( M-N\frac{\lambda}{\lambda+\mu} \right) \left( \frac{z^{\nu-1}}{z^\nu+(\lambda+\mu)}
					- \frac{z^{\nu-1}}{z^\nu + 2(\lambda+\mu)} \right)  \notag \\
					& + 2 \lambda^2 N(N-1) \frac{1}{\lambda+\mu} \cdot \frac{z^{-1}}{z^\nu
					+2(\lambda+\mu)} \notag.
				\end{align}
				Equation \eqref{elastic} then implies that
				\begin{align}
					H^\nu&(t) \\
					= {} & M(M-1) E_{\nu,1} \left( -2(\lambda+\mu)t^\nu \right) \notag \\
					& +\frac{2\lambda (N-1)}{\lambda+\mu} \left( M-N\frac{\lambda}{\lambda+\mu} \right)
					\left( E_{\nu,1}(-(\lambda+\mu)t^\nu) - E_{\nu,1}(-2(\lambda+\mu)t^\nu) \right) \notag \\
					& + 2\lambda^2 N(N-1)\frac{1}{\lambda+\mu} t^\nu E_{\nu,\nu+1} \left( -2(\lambda+\mu)t^\nu
					\right). \notag
				\end{align}
				Considering that $t^\nu E_{\nu,\nu+1}(at^\nu) = a^{-1} (E_{\nu,1}(at^\nu)-1)$, we obtain
				\begin{align}
					H^\nu(t) = {} & M(M-1) E_{\nu,1} \left( -2(\lambda+\mu)t^\nu \right) \\
					& +	\frac{2\lambda(N-1)}{\lambda+\mu} \left( M-N\frac{\lambda}{\lambda+\mu} \right)
					E_{\nu,1} (-(\lambda+\mu)t^\nu) \notag \\
					& - \frac{2\lambda (N-1)}{\lambda+\mu} \left( M-N\frac{\lambda}{\lambda+\mu} \right)
					E_{\nu,1} (-2(\lambda+\mu)t^\nu) \notag \\
					& -\frac{\lambda^2}{(\lambda+\mu)^2} N(N-1) E_{\nu,1} (-2(\lambda+\mu)t^\nu)
					+ \frac{\lambda^2}{(\lambda+\mu)^2} N(N-1) \notag \\
					= {} & M(M-1) E_{\nu,1} \left( -2(\lambda+\mu)t^\nu \right) \notag \\
					& +
					\frac{2\lambda M(N-1)}{\lambda+\mu} E_{\nu,1}(-(\lambda+\mu)t^\nu) \notag \\
					& - \frac{2 \lambda^2 N(N-1)}{(\lambda+\mu)^2} E_{\nu,1} (-(\lambda+\mu)t^\nu) \notag \\
					& - \frac{2\lambda M(N-1)}{\lambda+\mu} E_{\nu,1} (-2(\lambda+\mu)t^\nu) \notag \\
					& + \frac{2\lambda^2 N(N-1)}{(\lambda+\mu)^2} E_{\nu,1}(-2(\lambda+\mu)t^\nu) \notag \\
					& -\frac{\lambda^2 N(N-1)}{(\lambda+\mu)^2} E_{\nu,1}(-2(\lambda+\mu)t^\nu)
					+\frac{\lambda^2 N(N-1)}{(\lambda+\mu)^2} \notag \\
					= {} & \frac{\lambda^2 N(N-1)}{(\lambda+\mu)^2} \notag \\
					& + E_{\nu,1}(-2(\lambda+\mu)t^\nu)
					\left( \frac{\lambda^2 N(N-1)}{(\lambda + \mu)^2} - \frac{2\lambda M(N-1)}{\lambda+\mu}
					+M(M-1) \right) \notag \\
					& - E_{\nu,1} (-(\lambda+\mu)t^\nu) \left( \frac{2 \lambda^2 N(N-1)}{(\lambda+\mu)^2}
					-\frac{2\lambda M (N-1)}{\lambda+\mu} \right). \notag
				\end{align}
				The variance can thus be written as
				\begin{align}
					\mathbb{V}\text{ar} & \, \mathcal{N}^\nu(t) \\
					= {} & H^\nu(t) + \mathbb{E}\, \mathcal{N}^\nu(t)
					- (\mathbb{E}\, \mathcal{N}^\nu(t))^2 \notag \\
					= {} & H^\nu(t) + \left( M-N\frac{\lambda}{\lambda+\mu} \right) E_{\nu,1} (-(\lambda+\mu)t^\nu)
					\notag \\
					& + N \frac{\lambda}{\lambda+\mu} - \left( M-\frac{\lambda}{\lambda+\mu} \right)^2
					\left( E_{\nu,1}(-(\lambda+\mu)t^\nu) \right)^2 \notag \\
					& - N^2 \frac{\lambda^2}{(\lambda+\mu)^2} - 2 \frac{N \lambda}{\lambda+\mu}
					\left( M-N \frac{\lambda}{\lambda+\mu} \right) E_{\nu,1}(-(\lambda+\mu)t^\nu) \notag \\
					= {} & \left( \frac{\lambda^2 N(N-1)}{(\lambda +\mu)^2} -\frac{2\lambda M(N-1)}{\lambda+\mu}
					+ M(M-1) \right) E_{\nu,1} (-2(\lambda+\mu)t^\nu) \notag \\
					& + \left( \frac{2\lambda^2 N}{(\lambda+\mu)^2} - \frac{\lambda}{\lambda+\mu}
					(N+2M) +M \right)
					E_{\nu,1} (-(\lambda+\mu)t^\nu) \notag \\
					& - \left( M-N\frac{\lambda}{\lambda+\mu} \right)^2 \left( E_{\nu,1}(-(\lambda+\mu)t^\nu) \right)^2
					+ \frac{N \lambda\mu}{(\lambda+\mu)^2}. \notag
				\end{align}
			\end{proof}
			\begin{flushright} \qed \end{flushright}
		\end{thm}
		
		Exploiting Theorem \ref{sub}, we derive the explicit expression of the extinction
		probability $p_0^\nu(t) = \text{Pr} \{ \mathcal{N}^\nu(t) = 0 | \mathcal{N}^\nu(0) = M \}$ below.
		
		\begin{thm}
			The extinction probability $p_0^\nu(t) = \text{Pr} \{ \mathcal{N}^\nu(t) = 0 | \mathcal{N}^\nu(0) = M \}$
			for a fractional binomial process $\mathcal{N}^\nu(t)$, $t \ge 0$ is
			\begin{align}
				\label{extinction}
				p_0^\nu(t) = {} & \left( \frac{\mu}{\lambda+\mu} \right)^N \sum_{r=0}^{N-M} \binom{N-M}{r} \left( 
				\frac{\lambda}{\mu} \right)^r \\
				& \times \sum_{h=0}^M \binom{M}{h} (-1)^h E_{\nu,1}(-(r+h)(\lambda+\mu)t^\nu). \notag
			\end{align}
			
			\begin{proof}
				It is known \citep{jakeman} that the generating function $Q(u,t) = \sum_{n=0}^N (1-u)^n p_n(t)$
				for the classical binomial process can be written as
				\begin{align}
					Q(u,t)
					= {} & \left[ 1-\left( 1-e^{-(\mu+\lambda)t} \right)\frac{\lambda}{\lambda+\mu}u \right]^{N-M} \\
					& \times \left[ 1-\left( \left(1-e^{-(\mu+\lambda)t}\right)\frac{\lambda}{\lambda+\mu} +
					e^{-(\mu+\lambda)t} \right)u \right]^M. \notag
				\end{align}
				This suggests that the extinction probability for the classical case can be written as
				\begin{align}
					p_0(t) = {} & \left[1-\frac{\lambda}{\lambda+\mu}
					+ \frac{\lambda}{\lambda+\mu}
					e^{-(\mu+\lambda)t}\right]^{N-M} \\
					& \times \left[ 1-\left(\frac{\lambda}{\lambda+\mu}
					-e^{-(\mu+\lambda)t} \frac{\lambda}{\lambda+\mu}+e^{-(\mu+\lambda)t}\right) \right]^M \notag \\
					= {} & \left[ \frac{\mu}{\lambda+\mu} + \frac{\lambda}{\lambda+\mu} e^{-(\mu+\lambda)t} \right]^{N-M}
					\left[ \frac{\mu}{\lambda+\mu} - \frac{\mu}{\lambda+\mu}e^{-(\mu+\lambda)t} \right]^M \notag \\
					= {} & \left( \frac{1}{\lambda+\mu} \right)^N \mu^M \left( \mu+\lambda e^{-(\lambda+\mu)t} \right)^{N-M}
					\left( 1-e^{-(\lambda+\mu)t} \right)^M \notag \\
					= {} & \left( \frac{\mu}{\lambda+\mu} \right)^N \left( 1+\frac{\lambda}{\mu}
					e^{-(\lambda+\mu)t} \right)^{N-M} \left( 1-e^{-(\lambda+\mu)t} \right)^M \notag \\
					= {} & \left( \frac{\mu}{\lambda+\mu} \right)^N \sum_{r=0}^{N-M} \binom{N-M}{r} \left( 
					\frac{\lambda}{\mu} \right)^r e^{-r(\lambda+\mu)t} \notag \\
					& \times \sum_{h=0}^M
					\binom{M}{h} (-1)^h e^{-h(\lambda+\mu)t}
					\notag \\
					= {} & \left( \frac{\mu}{\lambda+\mu} \right)^N \sum_{r=0}^{N-M} \binom{N-M}{r} \left( 
					\frac{\lambda}{\mu} \right)^r \sum_{h=0}^M \binom{M}{h} (-1)^h e^{-(r+h)(\lambda+\mu)t}. \notag
				\end{align}
				Using Theorem \ref{sub}, we now obtain
				\begin{align}
					p_0^\nu(t) = {} & \int_0^\infty p_0(s) h(s,t) \mathrm ds \\
					= {} & \left( \frac{\mu}{\lambda+\mu} \right)^N
					\sum_{r=0}^{N-M} \binom{N-M}{r}\left( \frac{\lambda}{\mu} \right)^r \notag \\
					& \times \sum_{h=0}^M \binom{M}{h}
					(-1)^h \int_0^\infty e^{-(r+h)(\lambda+\mu)s} q(s,t) \mathrm ds \notag \\
					= {} & \left( \frac{\mu}{\lambda+\mu} \right)^N \sum_{r=0}^{N-M} \binom{N-M}{r} \left( 
					\frac{\lambda}{\mu} \right)^r \notag \\
					& \times \sum_{h=0}^M \binom{M}{h} (-1)^h
					E_{\nu,1}(-(r+h)(\lambda+\mu)t^\nu). \notag
				\end{align}
			\end{proof}
			\begin{flushright} \qed \end{flushright}
		\end{thm}

		\begin{thm}
			The state probabilities $p_n^\nu(t) = \text{Pr} \{ \mathcal{N}^\nu(t) = n | \mathcal{N}^\nu(0)
			= M \}$, $\lambda >0$, $\mu > 0$, have the following form:
			\begin{align}
				\label{stato}
				p_n^\nu(t) =
				\begin{cases}
					\sum_{r=0}^n g_{n,r}^\nu(t), & 0 \leq n < \min(M,N-M), \\
					\sum_{r=0}^{N-M} g_{n,r}^\nu(t), & N-M \leq n < M, \: M > N-M, \\
					\sum_{r=n-M}^n g_{n,r}^\nu(t), & M \leq n < N-M, \: M < N-M, \\
					\sum_{r=n-M}^{N-M} g_{n,r}^\nu(t), & \max(M,N-M) \leq n \leq N,
				\end{cases}
			\end{align}
			and
			\begin{align}
				p_n^\nu(t) =
				\begin{cases}
					\sum_{r=0}^n g_{n,r}^\nu(t), & 0 \leq n< M, \\
					\sum_{r=n-M}^M g_{n,r}^\nu(t), & M \leq n \leq N,
				\end{cases}			
			\end{align}
			when $N-M=M$, and where
			\begin{align}
				g_{n,r}^\nu(t)
				= {} & \left( \frac{\mu}{\lambda+\mu} \right)^N \binom{N-M}{r} \binom{M}{n-r} \\
				& \times \sum_{m_1=0}^r \binom{r}{m_1} (-1)^{m_1}
				\sum_{m_2=0}^{N-M-r} \binom{N-M-r}{m_2} \left( \frac{\lambda}{\mu} \right)^{m_2} \notag \\
				& \times \sum_{m_3=0}^{n-r} \binom{n-r}{m_3} \left( \frac{\lambda}{\mu} \right)^{n-m_3}
				\sum_{m_4=0}^{M-n+r} \binom{M-n+r}{m_4} (-1)^{m_4} \notag \\
				& \times E_{\nu,1} \left(
				-(m_1+m_2+m_3+m_4)(\mu+\lambda)t^\nu \right). \notag
			\end{align}
			
			\begin{proof}
				We start by rewriting  the probability generating function of the classical binomial
				process  as
				\begin{align}
					Q(u,t)& \\
					= {} & \left[ 1-\left( 1-e^{-(\mu+\lambda)t} \right)\frac{\lambda}{\lambda+\mu}u \right]^{N-M}\notag\\
					& \times \left[ 1-\left( \left(1-e^{-(\mu+\lambda)t}\right)\frac{\lambda}{\lambda+\mu} +
					e^{-(\mu+\lambda)t} \right)u \right]^M \notag \\
					= {} & \left[ (1-u) \left( \frac{\lambda}{\lambda+\mu} - \frac{\lambda}{\lambda+\mu}
					e^{-(\mu+\lambda)t}	\right) + \frac{\mu}{\lambda+\mu} + \frac{\lambda}{\lambda+\mu}
					e^{-(\mu+\lambda)t} \right]^{N-M} \notag \\
					& \times \left[ \frac{\lambda}{\lambda+\mu}(1-u) + \frac{\mu}{\lambda+\mu} + \frac{\mu}{\lambda+\mu}
					e^{-(\mu+\lambda)t} \right. \notag \\
					& \left. -\frac{\mu}{\lambda+\mu} e^{-(\mu+\lambda)t}
					-u \frac{\mu}{\lambda+\mu} e^{-(\mu+\lambda)t} \right]^M \notag \\
					= {} & \left( \frac{\lambda}{\lambda+\mu} \right)^N
					\left[ (1-u) \left( 1-e^{-(\mu+\lambda)t} \right) + \frac{\mu}{\lambda}
					+e^{-(\mu+\lambda)t} \right]^{N-M} \notag \\
					& \times \left[ (1-u) \left( 1+\frac{\mu}{\lambda}e^{-(\mu+\lambda)t} \right)
					+ \frac{\mu}{\lambda} - \frac{\mu}{\lambda} e^{-(\mu+\lambda)t} \right]^M \notag \\
					= {} & \left( \frac{\lambda}{\lambda+\mu} \right)^N
					\sum_{r=0}^{N-M} \binom{N-M}{r} (1-u)^r \notag \\
					& \times \left( 1-e^{-(\mu+\lambda)t} \right)^r
					\left( \frac{\mu}{\lambda} + e^{-(\mu+\lambda)t} \right)^{N-M-r} \notag \\
					& \times \sum_{h=0}^M \binom{M}{h} (1-u)^h \left( 1+\frac{\mu}{\lambda}
					e^{-(\mu+\lambda)t} \right)^h \left( \frac{\mu}{\lambda}-\frac{\mu}{\lambda}
					e^{-(\mu+\lambda)t} \right)^{M-h} \notag \\
					= {} & \left( \frac{\mu}{\lambda+\mu} \right)^N \sum_{r=0}^{N-M} \sum_{j=r}^{M+r}
					(1-u)^j \binom{N-M}{r} \binom{M}{j-r} \left( \frac{\lambda}{\mu}
					-\frac{\lambda}{\mu} e^{-(\mu+\lambda)t} \right)^r \notag \\
					& \times \left( 1+\frac{\lambda}{\mu} e^{-(\mu+\lambda)t} \right)^{N-M-r}
					\left( \frac{\lambda}{\mu} + e^{-(\mu+\lambda)t} \right)^{j-r} \notag \\
					& \times \left( 1-e^{-(\mu+\lambda)t} \right)^{M-j+r}. \notag
				\end{align}
			Letting
				\begin{align}
					g_{j,r}(t) = {} & \left( \frac{\mu}{\lambda+\mu} \right)^N
					\binom{N-M}{r} \binom{M}{j-r} \left( \frac{\lambda}{\mu}
					-\frac{\lambda}{\mu} e^{-(\mu+\lambda)t} \right)^r \\
					& \times \left( 1+\frac{\lambda}{\mu} e^{-(\mu+\lambda)t} \right)^{N-M-r}
					\left( \frac{\lambda}{\mu} + e^{-(\mu+\lambda)t} \right)^{j-r} \notag \\
					& \times \left( 1-e^{-(\mu+\lambda)t} \right)^{M-j+r}, \notag
				\end{align}
				we have
				\begin{align}
					Q(u,t) =
					\begin{cases}
						\sum_{j=0}^{M-1} (1-u)^j \sum_{r=0}^j g_{j,r}(t) \\
						\qquad + \sum_{j=M}^{N-M-1}
						(1-u)^j \sum_{r=j-M}^j g_{j,r}(t) \\
						\qquad + \sum_{j=N-M}^N (1-u)^j \sum_{j-M}^{N-M} g_{j,r}(t),
						& M < N-M, \\
						\sum_{j=0}^{N-M-1} (1-u)^j \sum_{r=0}^j g_{j,r}(t) \\
						\qquad + \sum_{j=N-M}^{M-1}
						(1-u)^j \sum_{r=0}^{N-M} g_{j,r}(t)\\
						\qquad + \sum_{j=M}^N (1-u)^j \sum_{r=j-M}^{N-M} g_{j,r}(t),
						& M > N-M, \\
						\sum_{j=0}^{M-1} (1-u)^j \sum_{r=0}^j g_{j,r}(t) \\
						\qquad + \sum_{j=M}^N (1-u)^j \sum_{r=j-M}^{M} g_{j,r}(t), & N-M=M.
					\end{cases}						
				\end{align}
				The classical state probabilities therefore read
				\begin{align}
					p_n(t) =
					\begin{cases}
						\sum_{r=0}^n g_{n,r}(t), & 0 \leq n < \min(M,N-M), \\
						\sum_{r=0}^{N-M} g_{n,r}(t), & N-M \leq n < M, \: M > N-M, \\
						\sum_{r=n-M}^n g_{n,r}(t), & M \leq n < N-M, \: M < N-M, \\
						\sum_{r=n-M}^{N-M} g_{n,r}(t), & \max(M,N-M) \leq n \leq N,
					\end{cases}
				\end{align}
				which reduce to
				\begin{align}
					p_n(t) =
					\begin{cases}
						\sum_{r=0}^n g_{n,r}(t), & 0 \leq n< M, \\
						\sum_{r=n-M}^M g_{n,r}(t), & M \leq n \leq N,
					\end{cases}			
				\end{align}
				when $N-M=M$.
				
				Exploiting Theorem \ref{sub}, we can derive the state probabilities for the
				fractional binomial process $\mathcal{N}^\nu(t)$, $t \ge 0$, as
				\begin{align}
					p_n^\nu(t) & = \int_0^\infty p_n(s) h(s,t) \mathrm ds \\
					& =
					\begin{cases}
						\sum_{r=0}^n g_{n,r}^\nu(t), & 0 \leq n < \min(M,N-M), \\
						\sum_{r=0}^{N-M} g_{n,r}^\nu(t), & N-M \leq n < M, \: M > N-M, \\
						\sum_{r=n-M}^n g_{n,r}^\nu(t), & M \leq n < N-M, \: M < N-M, \\
						\sum_{r=n-M}^{N-M} g_{n,r}^\nu(t), & \max(M,N-M) \leq n \leq N,
					\end{cases} \notag
				\end{align}
				or
				\begin{align}
					p_n^\nu(t) = \int_0^\infty p_n(s) h(s,t) \mathrm ds =
					\begin{cases}
						\sum_{r=0}^n g_{n,r}^\nu(t), & 0 \leq n< M, \\
						\sum_{r=n-M}^M g_{n,r}^\nu(t), & M \leq n \leq N,
					\end{cases}						
				\end{align}				
				for $N-M = M$. Note that
				\begin{align}
					g_{n,r}^\nu(t) = {} & \int_0^\infty g_{n,r}(s) h(s,t) \mathrm ds \\
					= {} & \int_0^\infty \left[ \left( \frac{\mu}{\lambda+\mu} \right)^N
					\binom{N-M}{r} \binom{M}{n-r} \left( \frac{\lambda}{\mu}
					-\frac{\lambda}{\mu} e^{-(\mu+\lambda)s} \right)^r \right. \notag \\
					& \left. \times \left( 1+\frac{\lambda}{\mu} e^{-(\mu+\lambda)s} \right)^{N-M-r}
					\left( \frac{\lambda}{\mu} + e^{-(\mu+\lambda)s} \right)^{n-r} \right. \notag \\
					& \times \left. \left( 1-e^{-(\mu+\lambda)s} \right)^{M-n+r} \right] h(s,t) \mathrm ds \notag \\
					= {} & \left( \frac{\mu}{\lambda+\mu} \right)^N \binom{N-M}{r} \binom{M}{n-r}
					\sum_{m_1=0}^r \binom{r}{m_1} (-1)^{m_1} \notag \\
					& \times \sum_{m_2=0}^{N-M-r} \binom{N-M-r}{m_2} \left( \frac{\lambda}{\mu} \right)^{m_2} \notag \\
					& \times \sum_{m_3=0}^{n-r} \binom{n-r}{m_3} \left( \frac{\lambda}{\mu} \right)^{n-m_3}
					\sum_{m_4=0}^{M-n+r} \binom{M-n+r}{m_4} (-1)^{m_4} \notag \\
					& \times \int_0^\infty e^{
					-(m_1+m_2+m_3+m_4)(\mu+\lambda)s} h(s,t) \mathrm ds \notag \\
					= {} & \left( \frac{\mu}{\lambda+\mu} \right)^N \binom{N-M}{r} \binom{M}{n-r}
					\sum_{m_1=0}^r \binom{r}{m_1} (-1)^{m_1} \notag \\
					& \times \sum_{m_2=0}^{N-M-r} \binom{N-M-r}{m_2} \left( \frac{\lambda}{\mu} \right)^{m_2} \notag \\
					& \times \sum_{m_3=0}^{n-r} \binom{n-r}{m_3} \left( \frac{\lambda}{\mu} \right)^{n-m_3}
					\sum_{m_4=0}^{M-n+r} \binom{M-n+r}{m_4} (-1)^{m_4} \notag \\
					& \times E_{\nu,1} \left(
					-(m_1+m_2+m_3+m_4)(\mu+\lambda)t^\nu \right). \notag
				\end{align}
				This concludes the proof.
			\end{proof}
			\begin{flushright} \qed \end{flushright}
		\end{thm}
		
		\begin{remark}
			From \eqref{stato},  we retrieve the extinction probability \eqref{extinction} when $n=0$.
		\end{remark}
		
		\begin{remark}
			As $t \rightarrow \infty$, the population in a fractional binomial
			process obeys a binomial distribution, i.e., 
			\begin{align}
				& \lim_{t \rightarrow \infty} Q^\nu(u,t) \\
				& = \sum_{r=0}^{N-M} \sum_{j=r}^{M+r} (1-u)^j \left( \frac{\mu}{\lambda+\mu} \right)^N
				\binom{N-M}{r} \binom{M}{j-r} \left( \frac{\lambda}{\mu} \right)^j \notag \\
				& = \sum_{r=0}^{N-M} \sum_{h=0}^M (1-u)^{h+r} \left( \frac{\mu}{\lambda+\mu} \right)^N
				\binom{N-M}{r} \binom{M}{h} \left( \frac{\lambda}{\mu} \right)^{h+r} \notag \\
				& = \left( \frac{\mu}{\lambda+\mu} \right)^N \left( 1+(1-u)\frac{\lambda}{\mu} \right)^{N-M}
				\left( 1+(1-u)\frac{\lambda}{\mu} \right)^M \notag \\
				& = \left( \frac{\mu}{\lambda+\mu} + \frac{\lambda}{\lambda+\mu}(1-u) \right)^N \notag \\
				& = \left( 1- \frac{\lambda}{\lambda+\mu}u \right)^N, \notag
			\end{align}
			is the probability generating function of a binomial random variable of parameter $\lambda/(\lambda+\mu)$
			and must be compared with equation $\mathrm{(9)}$ of \citet{jakeman}.
\end{remark}
Note also that as $t \rightarrow \infty$, indeed we observe (from \eqref{mea} and \eqref{va}) that
			\begin{align}
				\mathbb{E}\, \mathcal{N}^\nu(t) \longrightarrow N \frac{\lambda}{\lambda+\mu}
			\end{align}
and 
		\begin{align}
				\mathbb{V}ar\, \mathcal{N}^\nu(t) \longrightarrow N \frac{\lambda}{\lambda+\mu}\left(1- \frac{\lambda}{\lambda+\mu} \right),
			\end{align}
which are the mean and variance of the binomial equilibrium process. This suggests that the fractional generalization still preserves the binomial limit. 
		
	\section{Related fractional stochastic processes}
	
		In this section, we focus our attention to two pure branching processes which are in fact sub-models of the more general fractional binomial process described in Section \ref{se}. These are
		the fractional linear pure death process and the saturable fractional pure birth process. More specifically,	these processes can be directly obtained from the fractional binomial process by letting $\mu=0$ and 		$\lambda=0$, respectively. The main motivation underlining the analysis of these specific cases
		is that they are widely used in practice particularly in modeling  evolving populations 	in interacting environment possibly causing extintion or saturation.
		Our discussion on the fractional linear pure death process complements that of \citet{sakhno}'s. Instead, we analyze the saturable fractional pure birth
		process in more detail.
	
\begin{figure}[h!t!b!p!]
			\centering
			\includegraphics[height=2.5in, width=2.2in]{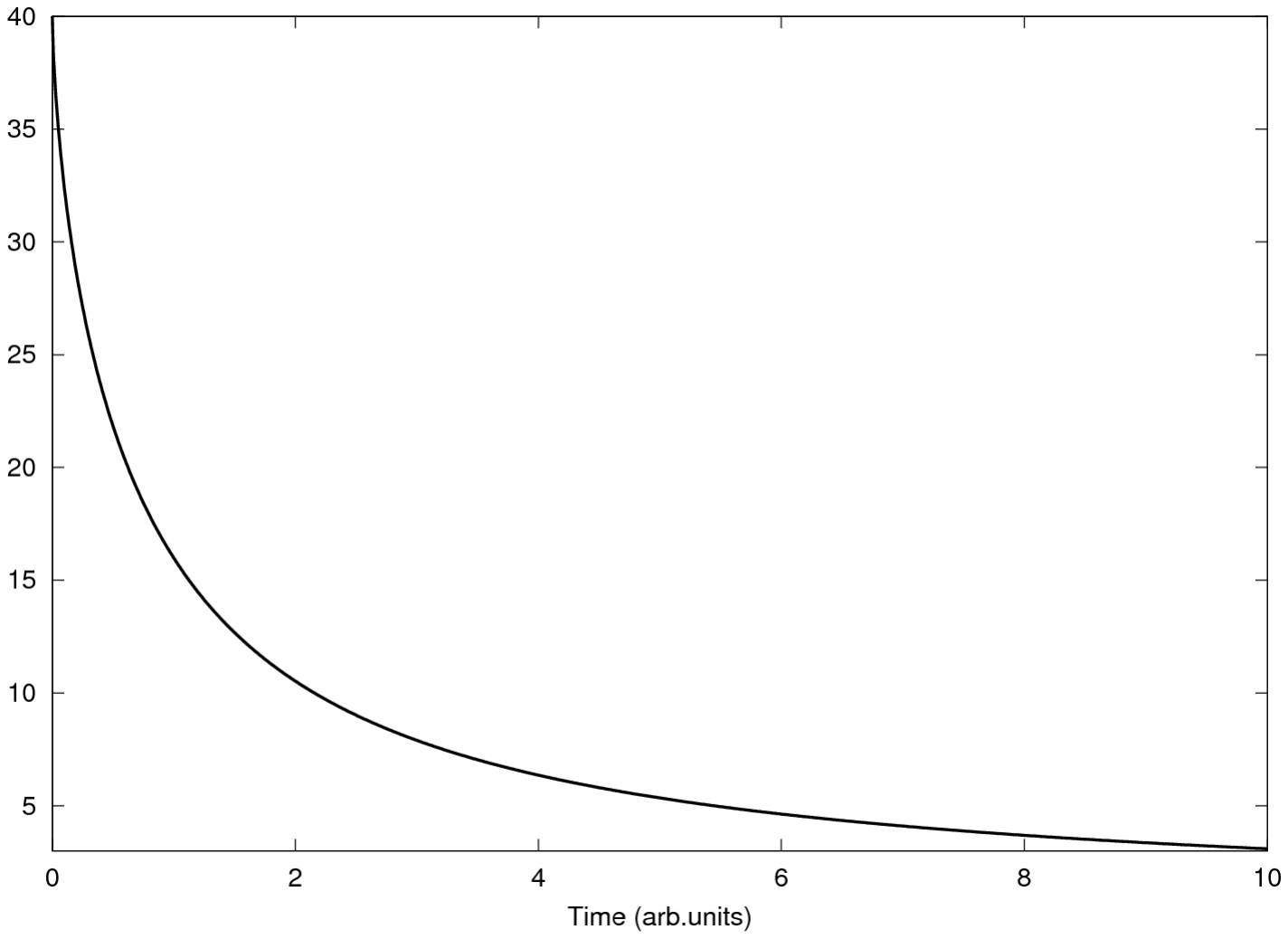}
			\includegraphics[height=2.5in, width=2.2in]{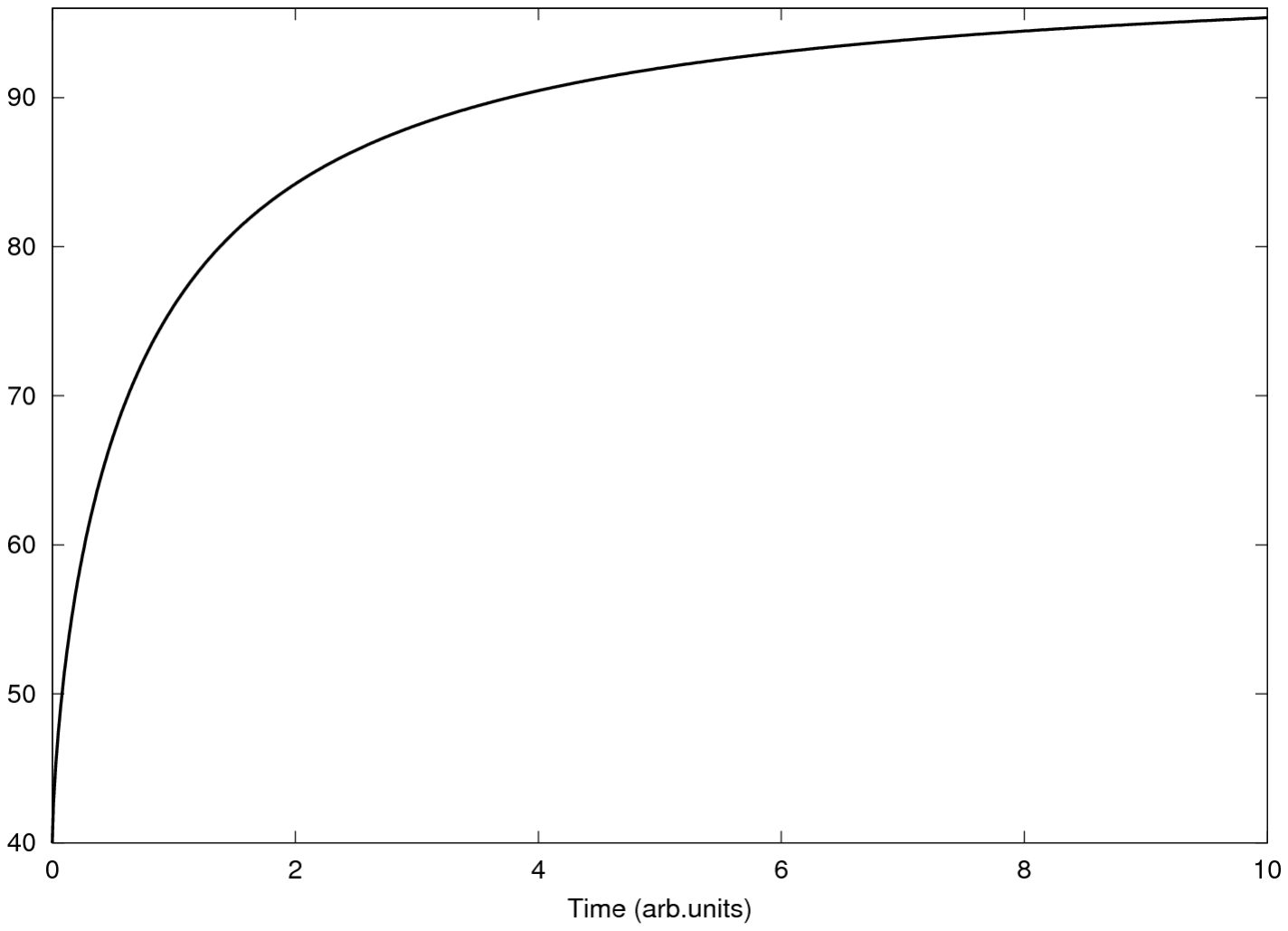}
			\caption{\label{bfig}The mean value of the fractional binomial process $\mathbb{E}\,
				\mathcal{N}^\nu(t)$ in the two different cases of pure death (left, $(\lambda,\mu)=(0,1)$)
				and pure birth (right, $(\lambda,\mu)=(1,0)$).
				For both cases we have $N=100$, $M=40$, $\nu=0.7$.}
		\end{figure}

		When $\lambda=0$, we obtain the mean value \eqref{mea} of the fractional linear pure
		death process $\mathcal{N}^\nu_d(t)$, $t \ge 0$, $\nu \in (0,1]$, (see \citet{sakhno}):
		\begin{align}
			\mathbb{E}\,\mathcal{N}^\nu_d(t) = M E_{\nu,1}(-\mu t^\nu), \qquad t \ge 0. 
		\end{align}
		Figure \ref{bfig} shows the mean value of the fractional linear pure death process (left)
		for specific values of the parameters $\mu$ and $\nu$.
		The variance can be  easily determined using \eqref{va} as
		\begin{align}
			\mathbb{V}\text{ar}\,\mathcal{N}^\nu_d(t) = {} & M(M-1) E_{\nu,1} \left( -2\mu t^\nu \right)
			+ M E_{\nu,1}\left( -\mu t^\nu \right) \\
			& - M^2 \left[ E_{\nu,1}\left( - \mu
			t^\nu \right) \right]^2, \qquad t \ge 0, \: \nu \in (0,1], \notag
		\end{align}
		and this reduces for $\nu = 1$ to the variance of the classical process (see \citet{bailey}, page 91, formula (8.32)):
		\begin{align}
			\mathbb{V}\text{ar}\,\mathcal{N}^1_d(t) = M e^{-\mu t}\left( 1-e^{-\mu t} \right).
		\end{align}
		Furthermore, the extinction probability of the fractional linear pure death
		process $\mathcal{N}^\nu_d(t)$, $t \ge 0$ (see \citet{sakhno}, page 73, formula (2.1)),
		\begin{align}
			\text{Pr} \{ \mathcal{N}^\nu_d(t) = 0 | \mathcal{N}^\nu_d(0) = M \}
			= \sum_{h=0}^M \binom{M}{h} (-1)^h E_{\nu,1} \left( -h \mu t^\nu \right), \qquad t \ge 0,
		\end{align}
		can also be derived directly from \eqref{extinction}.
		
		For the saturable fractional pure birth process $\mathcal{N}^\nu_b(t)$, 
		the mean value reduces to
		\begin{align}
			\label{satur}
			\mathbb{E}\,\mathcal{N}^\nu_b(t) = N - \left( N-M \right)E_{\nu,1}
			\left( -\lambda t^\nu \right).
		\end{align}
		Figure \ref{bfig} shows the expected value of the saturable fractional pure birth process (right)
		determined for specific values of the parameters $\lambda$ and $\nu$. The variance instead remains rather complicated and can be written by specialising \eqref{va} as
		\begin{align}
			\mathbb{V}\text{ar}\,\mathcal{N}^\nu_b(t) = {} &
			\left[ M(M-1)-N(N-1) \right] E_{\nu,1}\left( -2\lambda t^\nu \right) \\
			& - (N-M)(4N-1) E_{\nu,1}\left( -\lambda t^\nu \right) - (M-N)^2 \left[ E_{\nu,1}\left( -\lambda t^\nu
			\right) \right]^2. \notag
		\end{align}
As $t\rightarrow \infty$, 
		\begin{equation}
					\mathbb{E}\,\mathcal{N}^\nu_b(t) \longrightarrow   N 
		\end{equation}
		and
	\begin{equation}
					\mathbb{V}\text{ar}\,\mathcal{N}^\nu_b(t) \longrightarrow  0
		\end{equation}		
as expected.

We now determine the state probabilities $p_{n,b}^\nu(t) = \Pr \{ \mathcal{N}^\nu_b(t)=n
		| \mathcal{N}^\nu_b(0) = M \}$.
		When $\mu = 0$, the state probabilities can be derived from those of a
		nonlinear fractional pure birth process of \citet{pol},  and is given as
		\begin{equation}
			p_{n,b}^\nu(t) =
			\begin{cases}
				\prod_{j=M}^{n-1} \lambda_j \sum_{m=M}^n \frac{1}{\prod_{l=M,l \neq m}^n(\lambda_l-\lambda_m)}
				E_{\nu,1}\left( -\lambda_m t^\nu \right), & M < n \leq N ,\\
			E_{\nu,1} \left( -\lambda_M t^\nu \right), & n=M.
			\end{cases}
	\end{equation}
	Substituting the rates $\lambda_j= \lambda(N-j)$, we obtain
		\begin{align}
			p_{n,b}^\nu&(t) \\
			= {} & \prod_{j=M}^{n-1} \lambda(N-j) \sum_{m=M}^n
			\frac{E_{\nu,1}\left( -\lambda (N-m) t^\nu \right)}{\prod_{l=M,l \neq m}^n(\lambda (N-l)-\lambda(N-m))}
			\notag \\
			= {} & \sum_{m=M}^n \frac{(N-M)(N-M-1)\dots (N-n+1)}{(m-M)(m-M-1)\dots (m-m+1)
			(m-m-1)\dots (m-n)} \notag \\
			& \times E_{\nu,1}\left( -\lambda (N-m) t^\nu \right) \notag \\
			= {} & \sum_{m=M}^n \frac{(N-M)!}{(N-n)!(m-M)!(n-m)!} (-1)^{n-m}
			E_{\nu,1}\left( -\lambda (N-m) t^\nu \right) \notag \\
			= {} & \binom{N-M}{N-n} \sum_{m=M}^n \binom{n-M}{m-M} (-1)^{n-m}
			E_{\nu,1}\left( -\lambda (N-m) t^\nu \right), \quad M \leq n \leq N. \notag
		\end{align}
		This and formula \eqref{satur} show	that the behaviour of the saturable fractional pure birth
		process is subtantially different from that of the fractional Yule process. Similarly, the inter-birth waiting time $T_j^\nu$, i.e.\ the random
		time separating the $j$th and $(j+1)$th birth, has law
		\begin{align}
			\Pr \{ T_j^\nu \in \mathrm ds \} = \lambda(N-j) s^{\nu-1} E_{\nu,\nu} \left( -\lambda (N-j)s^\nu \right)
			\mathrm ds, \qquad j \ge M, \: s \ge 0.
		\end{align}
	The figure below shows the sample paths of the saturable fractional (bottom) and classical (top) linear pure birth processes. Apparently, the proposed model naturally includes processes or populations that saturate faster than the classical  linear pure birth process. The figure also indicates that saturation of the fractional binomial process is faster due to the explosive growth/birth bursts  at small times and as $\nu \to 0$.  Note that the parameters of these related fractional point processes can be estimated using the procedures of  \citet{cap11b}.
		\begin{figure}[h!t!b!p!]
			\centering
			\includegraphics[height=3in, width=4.5in]{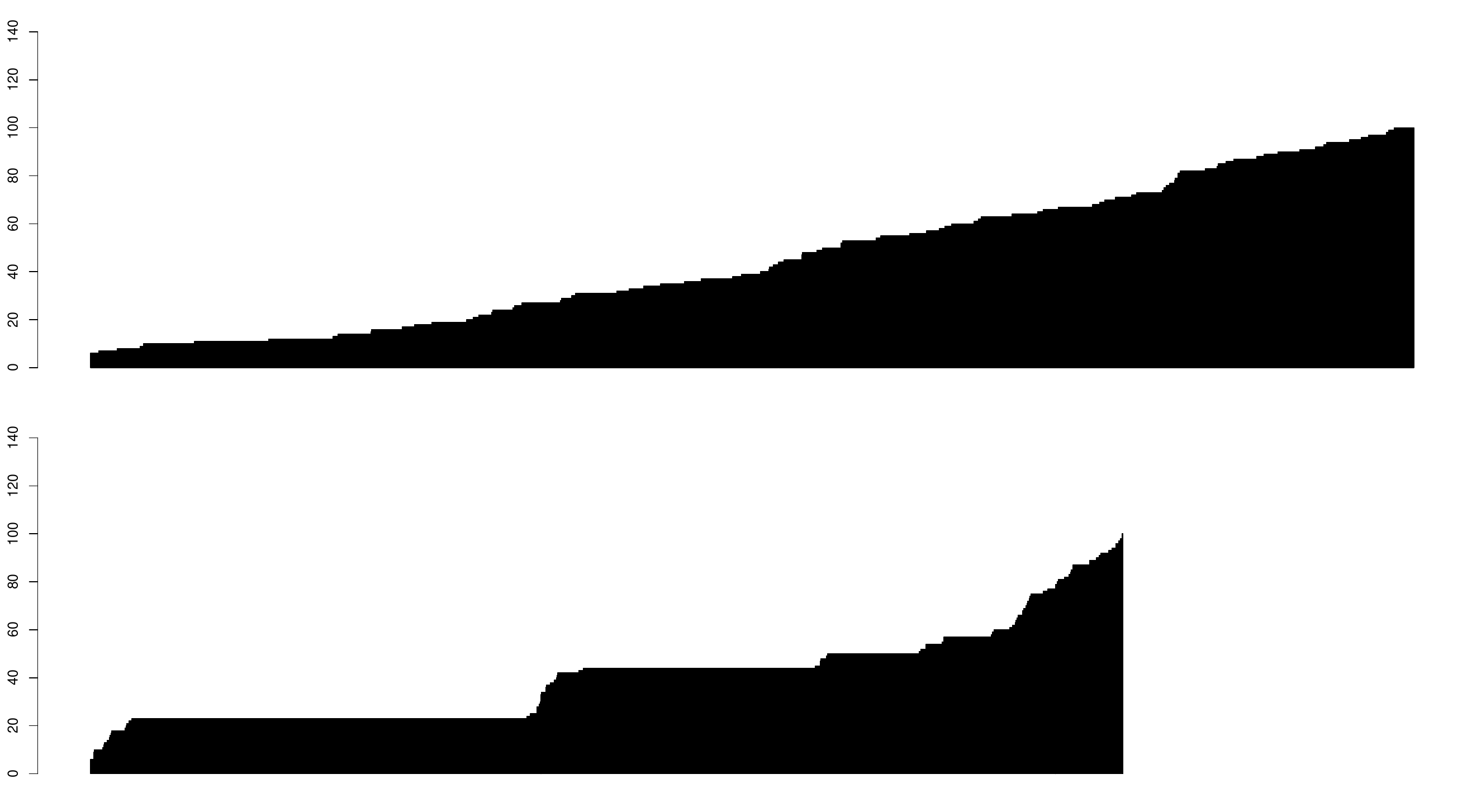}
			\caption{\label{bfig3}Sample trajectories of classical (top)  and fractional (bottom)  saturable  linear pure birth processes  using values $(N,M,\lambda, \nu)=(100,5,1,1)$  and $(N,M,\lambda, \nu)=(100,5,1,0.75)$, correspondingly.}
		\end{figure}

\section{Concluding Remarks}
    We have proposed a generalization of the binomial process using the techniques of fractional calculus. The fractional generalization In addition, more statistical properties of the fractional binomial process were derived. One interesting property of the fractional binomial process  was that it still preserved the binomial limit at large times while enlarging the class of models at  small and regular times that naturally include non-Markovian fluctuations with long memory. This potential made the proposed fractional binomial process appealing for real applications especially  to the quantum optics community.  New sub-models such as the saturable fractional pure birth process could also be automatically extracted from the proposed model. The generated sample trajectories of the saturable fractional linear pure birth process showed interesting features of the process such as the isolated bursts of the population growth particularly at small times. Overall, the fractional binomial process could be considered as a viable generalization of the classical binomial process.
    
    Although theoretical investigations have been done in the present paper, a number of issues are still left undone which could be considered as possible research extensions
of the current exploration. These may include: application of this model to rapidly saturable binomial processes, and the formalization of the parameter estimation procedures of the proposed model.


\end{document}